\documentclass{article}

\usepackage{lmodern} 
\usepackage{amsfonts} 
\usepackage{xfrac} 
\usepackage{siunitx} 
\usepackage{subcaption} 
\usepackage{multirow} 
\usepackage{amsmath}        
\usepackage{amsfonts}       
\usepackage{bm}             
\usepackage{mathcommand}    
\usepackage{mathtools}      
\usepackage{xparse}         

\renewcommand*{\vec}[1]{{\bm{#1}}}

\newcommand*{\ve}{\vec{e}}

\newcommand*{\vx}{\vec{x}}
\newcommand*{\vy}{\vec{y}}

\newcommand*{\vxi}{\vec{\xi}}

\newcommand*{\vnull}{\vec{0}}

\newcommand*{\mat}[1]{{\bm{#1}}}

\newcommand*{\mI}{\mat{I}}

\newcommand*{\mX}{\mat{X}}

\newcommand*{\op}[1]{{\mathcal{#1}}}

\newcommand*{\opY}{\op{Y}}

\newcommand*{\super}[2]{{#1^{#2}}}

\newcommand*{\transposesymbol}{\mathsf{T}}
\newcommand*{\hermitiansymbol}{\mathsf{H}}
\newcommand*{\adjointsymbol}{'}
\newcommand*{\inversesymbol}{-}

\newcommand*{\transpose}[1]{\super{#1}{\transposesymbol}}

\newcommand*{\hermitian}[1]{\super{#1}{\hermitiansymbol}}

\newcommand*{\T}[1]{\transpose{#1}}
\renewmathcommand{\H}[1]{\hermitian{#1}}

\newcommand*{\Tinv}[1]{\super{#1}{\inversesymbol\transposesymbol}}
\newcommand*{\adjinv}[1]{\super{#1}{\inversesymbol\adjointsymbol}}

\newcommand*{\invT}[1]{\Tinv}
\newcommand*{\invadj}[1]{\adjinv}
\newcommand*{\pinvT}[1]{\Tpseudoinv}
\newcommand*{\pinvadj}[1]{\adjpseudoinv}

\newcommand{\bsmallmat}[1]{
\begin{bsmallmatrix}
    #1
\end{bsmallmatrix}
}

\newcommand*{\dif}[2][1]{\mathrm{d}\ifstrequal{#1}{1}{}{^{#1}} #2}
\newcommand*{\pdif}[2][1]{\partial\ifstrequal{#1}{1}{}{^{#1}} #2}
\newcommand*{\diff}[2][1]{\mathrm{d} #2\ifstrequal{#1}{1}{}{^{#1}}}
\newcommand*{\pdiff}[2][1]{\partial #2\ifstrequal{#1}{1}{}{^{#1}}}

\newcommand*{\dx}{\diff{x}}

\newcommand*{\uniform}[2]{\mathcal{U}\left(#1, #2\right)}
\newcommand*{\normal}[2]{\mathcal{N}\left(#1, #2\right)}

\newcommand*{\probabilityof}[2][]{p\left(#2 \ifstrempty{#1}{}{|#1} \right)}
\newcommand*{\proposalof}[1]{q\left(#1\right)}
\newcommand*{\expectationof}[1]{\mathbb{E}\left[#1\right]}

\newcommand{\solution}{u}
\newcommand{\forcing}{f}
\newcommand{\diffusion}{\kappa}

\newcommand{\param}{\xi}
\newcommand{\params}{\vxi}
\newcommand{\eigenval}{\lambda}

\newcommand{\eigenfunc}{\phi}

\newcommand{\gradient}{\nabla}
\newcommand{\divergence}{\nabla \cdot}

\newcommand{\ofx}{(x)}
\newcommand{\ofvx}{(\vx)}
\newcommand{\domain}{\Omega}
\newcommand{\boundary}{\partial\Omega}
\newcommand{\boundarydir}{\Gamma_d}
\newcommand{\boundaryneu}{\Gamma_n}

\newcommand{\reference}[1]{{#1}^*}
\newcommand{\refsolution}{\reference{\solution}}

\newcommand{\refparams}{\reference{\params}}

\newcommand{\fem}[1]{{#1}^\text{h}}
\newcommand{\femsolution}{\fem{\solution}}
\newcommand{\vecsolution}{\vec{\solution}}
\newcommand{\vecforcing}{\vec{\forcing}}

\newcommand{\nodalerror}{\ve}
\newcommand{\fielderrornorm}{\|e\ofx\|}
\newcommand{\nodalerrornorm}{\|\ve\|}
\newcommand{\nodalerrorrel}{\zeta}
\newcommand{\intererrorrel}{\eta}

\newcommand{\obslocation}{X}
\newcommand{\obslocations}{\mX}

\newcommand{\perturbed}[1]{{\tilde{#1}}}
\newcommand{\pertsolution}{\perturbed{\solution}}
\newcommand{\pertfemsolution}{\fem{\perturbed{\solution}}}
\newcommand{\pertlocation}{\perturbed{\obslocation}}
\newcommand{\pertlocations}{\perturbed{\obslocations}}

\newcommand{\observations}{\vy}
\newcommand{\obsoperator}{\opY}
\newcommand{\obsnoise}{\sigma_e}

\newcommand{\prior}{\probabilityof{\params}}
\newcommand{\likelihood}{\probabilityof[\params]{\observations}}
\newcommand{\likelihoodu}{\probabilityof[\solution\ofx]{\observations}}
\newcommand{\posterior}{\probabilityof[\observations]{\params}}
\newcommand{\posteriorkappa}{\probabilityof[\observations]{\diffusion\ofx}}
\newcommand{\evidence}{\probabilityof{\observations}}

\newcommand{\nelem}{n}
\newcommand{\nmetropolis}{N}
\newcommand{\nmontecarlo}{M}

\newcommand{\meshsize}{h}

\newcommand{\energy}{E}
\newcommand{\refenergy}{\reference{\energy}}
\newcommand{\femenergy}{\fem{\energy}}
\newcommand{\pertfemenergy}{\fem{\perturbed{\energy}}}

\usepackage{xspace}

\newcommand{\etalia}{et al\futurelet\temp\addspaceandorperiod}

\newcommand{\idest}{i.e.,\ }
\newcommand{\ie}{\idest}
\newcommand{\exempligratia}{e.g.,\ }
\newcommand{\eg}{\exempligratia}
\newcommand{\etcetera}{etc\futurelet\temp\addspaceandorperiod}

\newcommand{\independentidenticallydistributed}{i.i.d\futurelet\temp\addspaceandorperiod}
\newcommand{\iid}{\independentidenticallydistributed}

\def\addspaceandorperiod{%
  \ifx.\temp{}
  \else%
    \ifcat a\temp\relax{.\ }
    \else%
      \ifcat A\temp\relax{.\ }
      \else%
        {.\nobreak\hspace{0pt}}
      \fi
    \fi
  \fi
}

\newcommand{\folder}{4-observations}

\usepackage[final]{ProbNum25}

\usepackage[capitalise,nameinlink]{cleveref}

\probnumtitle{Effects of Interpolation Error and Bias on the Random Mesh Finite Element Method for Inverse Problems}
\probnumauthors{%
    \name{Anne Poot} \affiliation{Delft University of Technology, Delft, The Netherlands} \and%
    \name{Iuri Rocha} \affiliation{Delft University of Technology, Delft, The Netherlands} \and%
    \name{Pierre Kerfriden} \affiliation{Mines Paris--PSL, \'Evry, France} \and%
    \name{Frans van der Meer} \affiliation{Delft University of Technology, Delft, The Netherlands}\and%
}

\probnumabstract{
    Bayesian inverse problems are an important application for probabilistic solvers of partial differential equations:
    when fully resolving numerical error is computationally infeasible, probabilistic solvers can be used to consistently model the error and propagate it to the posterior.
    In this work, the performance of the random mesh finite element method (RM-FEM) is investigated in a Bayesian inverse setting.
    We show how interpolation error negatively affects the RM-FEM posterior, and how these negative effects can be diminished.
    In scenarios where FEM is biased for a quantity of interest, we find that RM-FEM struggles to accurately model this bias.
}

\begin{document}

\section{Introduction}
\label{sec:introduction}
Within the field of probabilistic numerics, one key area of interest is the probabilistic modeling of partial differential equation (PDE) solver error.
In the past decade, many approaches have been proposed, which can be roughly divided into a Bayesian and a frequentist category.
In the Bayesian category, a (usually Gaussian process) prior distribution is assumed over the solution space, which is conditioned on the PDE's strong-form residuals \cite{raissi_machine_2017,cockayne_probabilistic_2017,yang_b-pinns_2021}.
Many classical solvers of PDEs can be reformulated in a Bayesian manner~\cite{pfortner_physics-informed_2023}, including the method of lines~\cite{kramer_probabilistic_2022}, the finite element method (FEM)~\cite{poot_bayesian_2024} and the finite volume method~\cite{weiland_scaling_2024}.
Methods in the frequentist category, on the other hand, typically rely on perturbation of the model in order to induce a distribution over model predictions to capture the model error.
Although these solver-perturbing methods are often focused on solving ordinary differential equations~\cite{lie_strong_2019,abdulle_random_2020}, they can be extended to PDEs as well \cite{conrad_statistical_2017}.

In this work, we focus on a particular perturbation-based PDE solver presented by Abdulle and Garegnani \cite{abdulle_probabilistic_2021}, which aims to model finite element discretization error by randomizing the node locations of the finite element mesh.
The method they present is appealing both from a theoretical and a practical point of view.
Theoretically, they are able to derive two error estimators from a second moment of the distribution over the gradients of the solution, by proving equivalence to the classic error estimator from \cite{babuska_analysis_1979}.
The practical appeal of the method lies in its ease of implementation:
because the distribution of $\solution\ofx$ only depends on randomized node locations, samples can easily be computed with off-the-shelf finite element software.

One driving motivation behind the development of probabilistic approaches to FEM is their potential application to (Bayesian) inverse problems.
The typical procedure to infer a set of parameters of the PDE from a set of observations involves Monte Carlo sampling, with each sample requiring an FEM solve.
To keep this approach computationally feasible, it might be necessary to resort to a coarse finite element discretization.
However, this will reduce the accuracy of both the FEM solution and the resulting parameter estimates.
The aim of taking a probabilistic approach to FEM is not necessarily to improve the accuracy of the posterior mean, but rather to incorporate the uncertainty due to discretization error in the posterior variance.
We say that a probabilistic formulation of FEM is performing well in an inverse setting if
1) its posterior mean is at least as close to the ground truth as the posterior obtained with traditional FEM,
and 2) its posterior variance accurately reflects the mismatch between the posterior mean and ground truth due to discretization error.
The aim of this work is to investigate the performance of RM-FEM in the context of a Bayesian inverse problem, focussing specifically on the effects of interpolation error and solver bias.

The paper is structured as follows:
in \cref{sec:model-definition}, we present the base model used throughout the paper and compare the results of FEM and RM-FEM in an inverse setting.
The influence of interpolation error and bias on the RM-FEM posterior are investigated in \cref{sec:interpolation-error,sec:bias}, respectively.
Finally, we present our conclusions in \cref{sec:conclusion}.

\section{Base model}
\label{sec:model-definition}
The analysis in this work is largely built on the work of Abdulle and Garegnani \cite[Section 4.2.1]{abdulle_probabilistic_2021}.
Before expanding on it, we first present their setup as a base model, and reproduce their results for both FEM and RM-FEM.

\subsection{Forward model}
\label{subsec:forward-model}
\begin{figure}
    \includegraphics[width=\columnwidth]{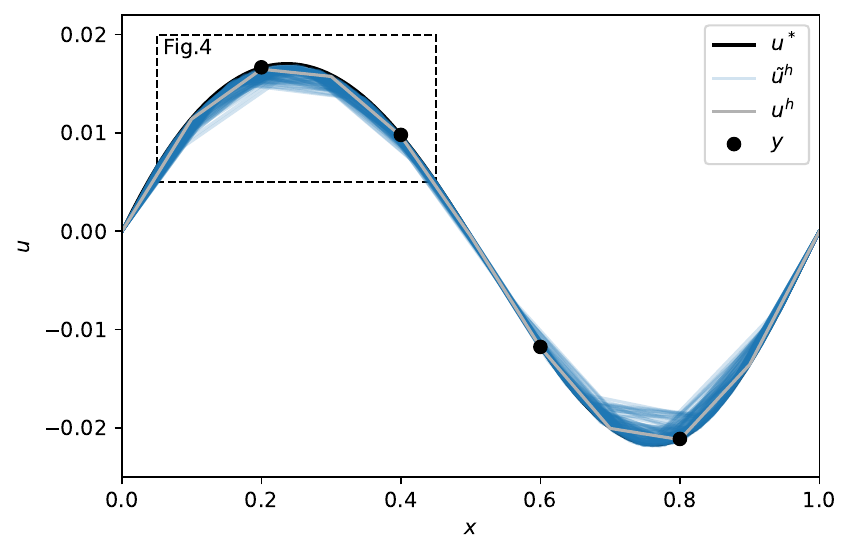}
    \caption{
        Comparison of the FEM solution $\femsolution\ofx$ and 100 samples from RM-FEM $\pertfemsolution\ofx$.
        The reference solution $\refsolution\ofx$ is shown in black, as well as the 4 observations $\observations$ that will be used in the inverse problem.
        All fields have been obtained using the reference parameters $\refparams$.
    }
    \label{fig:state0-plot}
\end{figure}

We consider the following 1D diffusion equation:
\begin{equation}
    \label{eq:strong-form}
    \begin{aligned}
        - \divergence \left(\diffusion\ofx \gradient \solution \ofx \right) &= \forcing \ofx & \forall x &\in \domain \\
        \solution \ofx &= 0 & \forall x &\in \boundary
    \end{aligned}
\end{equation}
on the domain $\domain = (0, 1)$.
The diffusion coefficient $\diffusion \ofx$ is parametrized by:
\begin{equation}
    \label{eq:diffusion}
    \diffusion\ofx = \exp\left\{\sum_{k=1}^4 \frac{\param_k}{\sqrt{\eigenval_k}} \eigenfunc_k\ofx\right\}
\end{equation}
where $\{\eigenval_k,\eigenfunc_k\ofx\}$ is the $k$th eigenpair of the differential operator $-\frac{\diff{}}{\diff[2]{x}}$, \ie $\eigenval_k = k^2 \pi^2$ and $\eigenfunc_k = \sqrt{2} \sin(k \pi x)$.
We take the forcing term to be $\forcing \ofx = \sin(2 \pi x)$ and set the parameter vector of our reference solution to $\params = \T{\bsmallmat{\param_1 & \param_2 & \param_3 & \param_4}} = \T{\bsmallmat{1 & 1 & \sfrac{1}{4} & \sfrac{1}{4}}}$.
The corresponding reference solution $\refsolution \ofx$ is shown in black in \cref{fig:state0-plot}.

In order to solve \cref{eq:strong-form}, we discretize the domain using $\nelem$ elements with linear shape functions.
A uniform mesh is used, so the element size $\meshsize$ is given simply by $\meshsize = \sfrac{1}{\nelem}$.
For the assembly of the stiffness matrix and force vector, a fourth-order Gauss-Legendre integration scheme is used.
In \cref{fig:state0-plot}, the finite element solution using $\nelem=10$ elements is shown in gray.

\subsection{Inverse model}
\label{sec:inverse-model}
We now formulate the inverse problem, aiming to infer $\diffusion\ofx$ from a set of observations $\observations$.
Our observation model is defined by:
\begin{align}
    \label{eq:observation}
    \observations &= \obsoperator[\solution\ofx] + \beta & \beta &\sim \normal{\vnull}{\sigma_e^2 \mI}
\end{align}
Here, the observation operator $\obsoperator$ takes a solution field $\solution \ofx$ as input, and returns the evaluation of that field at 4 equally spaced locations $\mX = \{\sfrac{i}{5}\}_{i=1}^4$.
The observations $\observations$ are generated by applying this operator to the reference solution $\refsolution \ofx$ and corrupting the result with \iid Gaussian noise ($\sigma_e = 10^{-5}$).
Compared to \cite{abdulle_probabilistic_2021}, the observation noise is reduced to make the effect of discretization error more visible:
for $\obsnoise=\num[print-zero-exponent=true, print-unity-mantissa=false]{1e-4}$, the effect of discretization error on the posterior is too small to be noticeable compared to the uncertainty due to observation noise.
We use \num{4} instead of \num{9} observation locations in order to better show the effect of excluding observation node locations from perturbation in \cref{subsec:fixed-observations}.

For our Bayesian inverse model, we assume a prior distribution over $\diffusion \ofx$ via the parametrization in \cref{eq:diffusion} by putting a prior on its parameters:
\begin{equation}
    \label{eq:prior}
    \probabilityof{\params} = \normal{\vnull}{\mI}
\end{equation}
The likelihood is defined according to \cref{eq:observation}, with $\sigma_e = 10^{-5}$.
Since $\solution\ofx$ is fully determined by $\params$, we have:
\begin{equation}
    \label{eq:likelihood-fem}
    \likelihood = \likelihoodu = \normal{\obsoperator\left[\solution\ofx\right]}{\sigma_e^2 \mI}
\end{equation}
Note that the model used to fit the data is almost exactly the same model as the one used to generate the data, with the notable difference that the data generation model uses $\reference{\nelem}=1000$ elements, whereas the fitted model uses a much coarser mesh with $\nelem=\{10, 20, 40\}$.
Although such similarity between the data generation model and Bayesian inversion model can be considered bad practice~\cite{wirgin_inverse_2004}, in our case it is a deliberate choice that helps distinguish the influence of discretization from other forms of model misspecification.

\subsection{Markov chain Monte Carlo}
\label{subsec:mcmc}
The posterior distribution is given by Bayes' theorem:
\begin{equation}
    \label{eq:posterior}
    \posterior = \frac{\prior \likelihood}{\evidence}
\end{equation}
In order to sample from $\probabilityof[\observations]{\params}$ we employ a random walk Metropolis algorithm~\cite{metropolis_equation_1953}.
The algorithm is started at $\params_0 = \vnull$, and the initial proposal distribution $\proposalof{\params}$ is set to be equal to the prior distribution $\probabilityof{\params}$.
We employ a burn-in period of $\num{10000}$ samples, during which the covariance of the proposal distribution is scaled isotropically to obtain a good acceptance ratio.
After the burn-in period, the covariance of the proposal distribution is fixed, the burn-in samples are discarded, and another $\nmetropolis = \num{10000}$ samples are drawn.
All results presented in this paper relate only to this second set of samples.

\subsection{RM-FEM}
\label{subsec:rmfem}
For our random mesh finite element implementation, we perturb the node locations according to \cite[Equation 2.4]{abdulle_probabilistic_2021}:
\begin{equation}
    \label{eq:perturbation}
    \pertlocation_i = \obslocation_i + h^p \alpha_i
\end{equation}
where $h = \sfrac{1}{\nelem}$, and $\alpha_i \sim \uniform{-\sfrac{1}{2}}{\sfrac{1}{2}}$.
We set $p=1$ in light of \cite[Theorem 2.9]{abdulle_probabilistic_2021}.
This distribution over node locations induces a distribution over FEM solutions $\pertsolution\ofx$, of which 100 samples are shown in blue in \cref{fig:state0-plot}.

In the random mesh setting, $\solution\ofx$ is now no longer fully determined by $\params$.
Therefore, we need to marginalize over all random meshes $\pertlocations$ to compute the likelihood:
\begin{equation}
    \label{eq:likelihood-rmfem}
    \likelihood = \int_\pertlocations \probabilityof[\params, \pertlocations]{\observations} \diff \pertlocations
\end{equation}
There are two approaches to perform this marginalization in an inverse setting, laid out in \cite{garegnani_sampling_2021}.
The first approach, Monte Carlo within Metropolis \cite{beaumont_estimation_2003}, replaces the integral in \cref{eq:likelihood-rmfem} by a sample approximation:
\begin{equation}
    \label{eq:likelihood-mcwm}
    \likelihood \approx \frac{1}{\nmontecarlo} \sum_{j=1}^\nmontecarlo \probabilityof[\params, \pertlocations_j]{\observations}
\end{equation}
In \cite{andrieu_pseudo-marginal_2009}, the resulting approximate posterior distribution is shown to converge asymptotically to the true posterior $\posterior$ as $\nmontecarlo$ grows large.
The second approach, Metropolis within Monte Carlo, performs the marginalization by approximating the posterior instead:
\begin{equation}
    \label{eq:posterior-mwmc}
    \posterior \approx \frac{1}{\nmontecarlo} \sum_{j=1}^\nmontecarlo \probabilityof[\observations, \pertlocations_j]{\params}
\end{equation}
This involves a full random walk Metropolis run on each perturbed mesh, and averaging the resulting $\nmontecarlo$ sets of samples obtain an ensemble approximation of $\posterior$.

Although the latter approach is easier to implement than the former, it tends to produce a rough multimodal approximation of $\posterior$ for the problem we are interested in here.
For this reason, we opt for the first approach, and approximate the likelihood $\likelihood$ according to \cref{eq:likelihood-mcwm}, with $\nmontecarlo=10$.

\subsection{Results}
\label{subsec:results}
\begin{figure*}
    \begin{subfigure}{0.3\textwidth}
        \includegraphics[width=\textwidth]{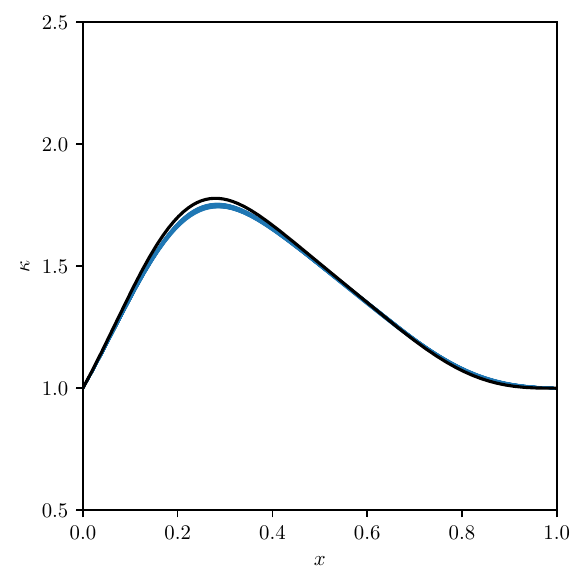}
        \caption{FEM, $\meshsize=\sfrac{1}{10}$}
        \label{subfig:kappa-plot-fem-10}
    \end{subfigure}
    \hfill
    \begin{subfigure}{0.3\textwidth}
        \includegraphics[width=\textwidth]{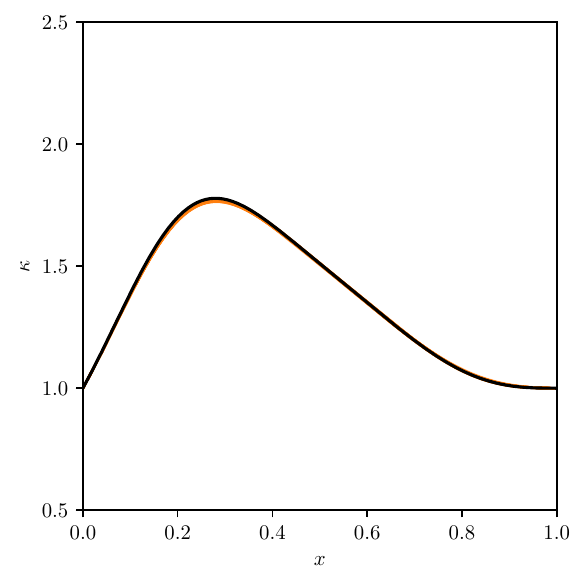}
        \caption{FEM, $\meshsize=\sfrac{1}{20}$}
        \label{subfig:kappa-plot-fem-20}
    \end{subfigure}
    \hfill
    \begin{subfigure}{0.3\textwidth}
        \includegraphics[width=\textwidth]{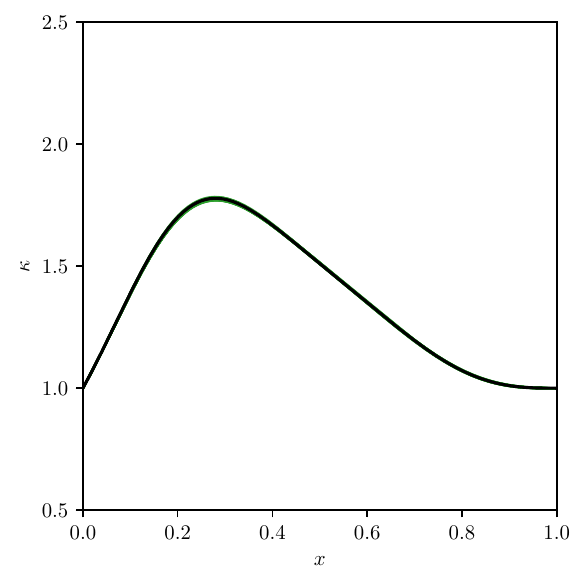}
        \caption{FEM, $\meshsize=\sfrac{1}{40}$}
        \label{subfig:kappa-plot-fem-40}
    \end{subfigure}

    \begin{subfigure}{0.3\textwidth}
        \includegraphics[width=\textwidth]{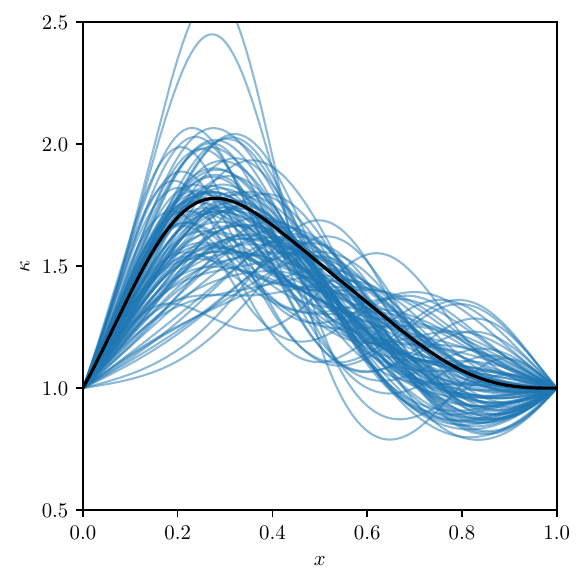}
        \caption{RM-FEM, $\meshsize=\sfrac{1}{10}$}
        \label{subfig:kappa-plot-rmfem-10}
    \end{subfigure}
    \hfill
    \begin{subfigure}{0.3\textwidth}
        \includegraphics[width=\textwidth]{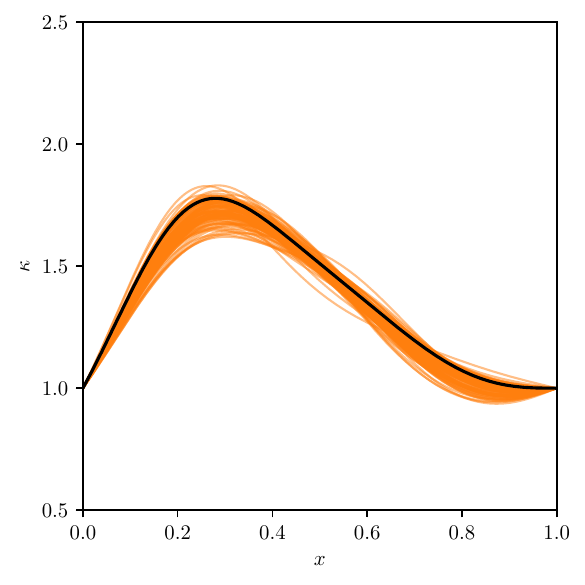}
        \caption{RM-FEM, $\meshsize=\sfrac{1}{20}$}
        \label{subfig:kappa-plot-rmfem-20}
    \end{subfigure}
    \hfill
    \begin{subfigure}{0.3\textwidth}
        \includegraphics[width=\textwidth]{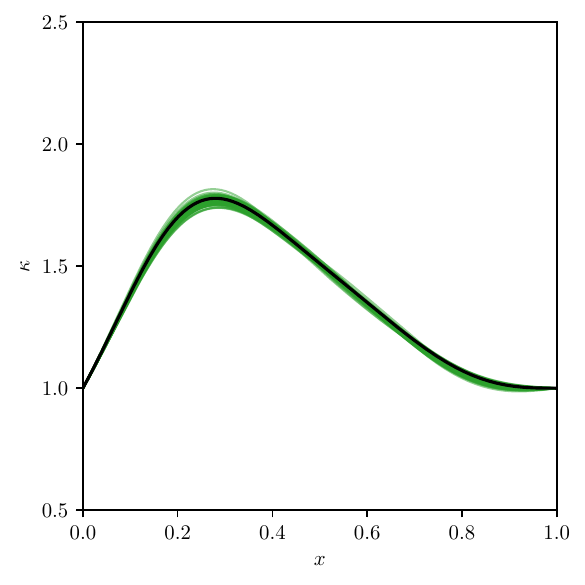}
        \caption{RM-FEM, $\meshsize=\sfrac{1}{40}$}
        \label{subfig:kappa-plot-rmfem-40}
    \end{subfigure}
    \caption{
        Reconstruction of Figure 10 from \cite{abdulle_probabilistic_2021}.
    }
    \label{fig:kappa-plots}
\end{figure*}

\begin{figure*}
    \begin{subfigure}{0.5\textwidth}
        \includegraphics[width=\textwidth]{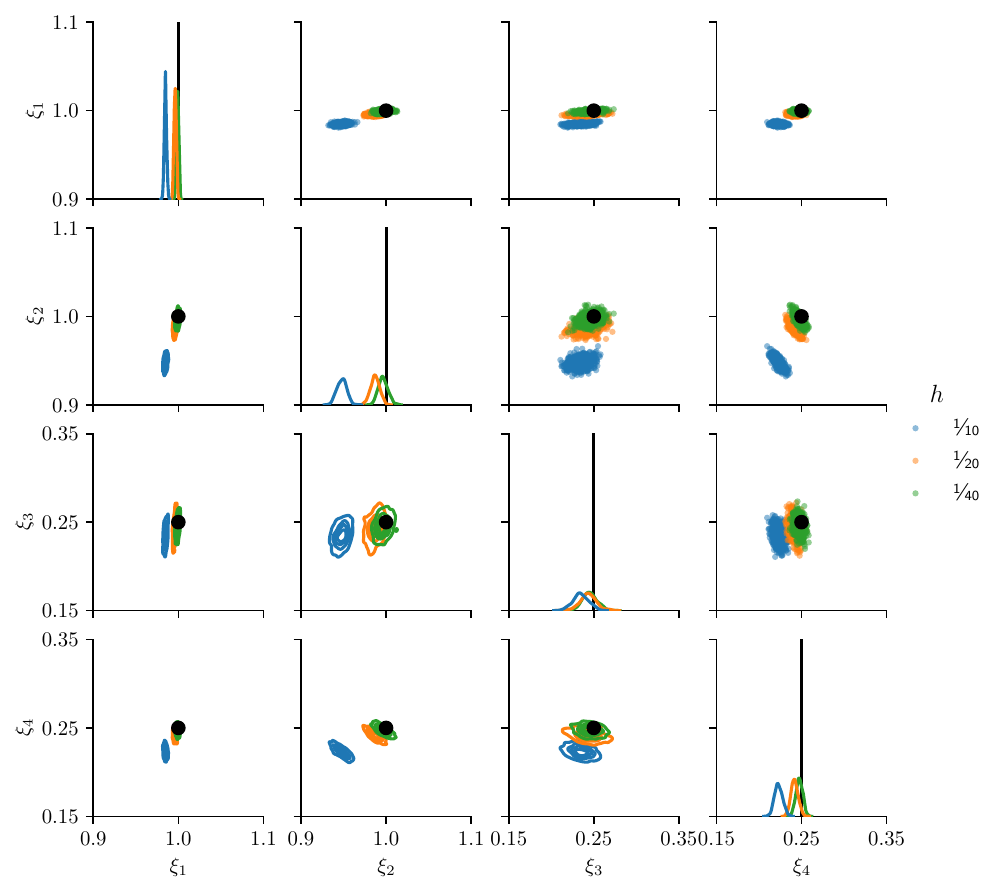}
        \caption{FEM}
        \label{subfig:pairgrid-plot-fem}
    \end{subfigure}
    \hfill
    \begin{subfigure}{0.5\textwidth}
        \includegraphics[width=\textwidth]{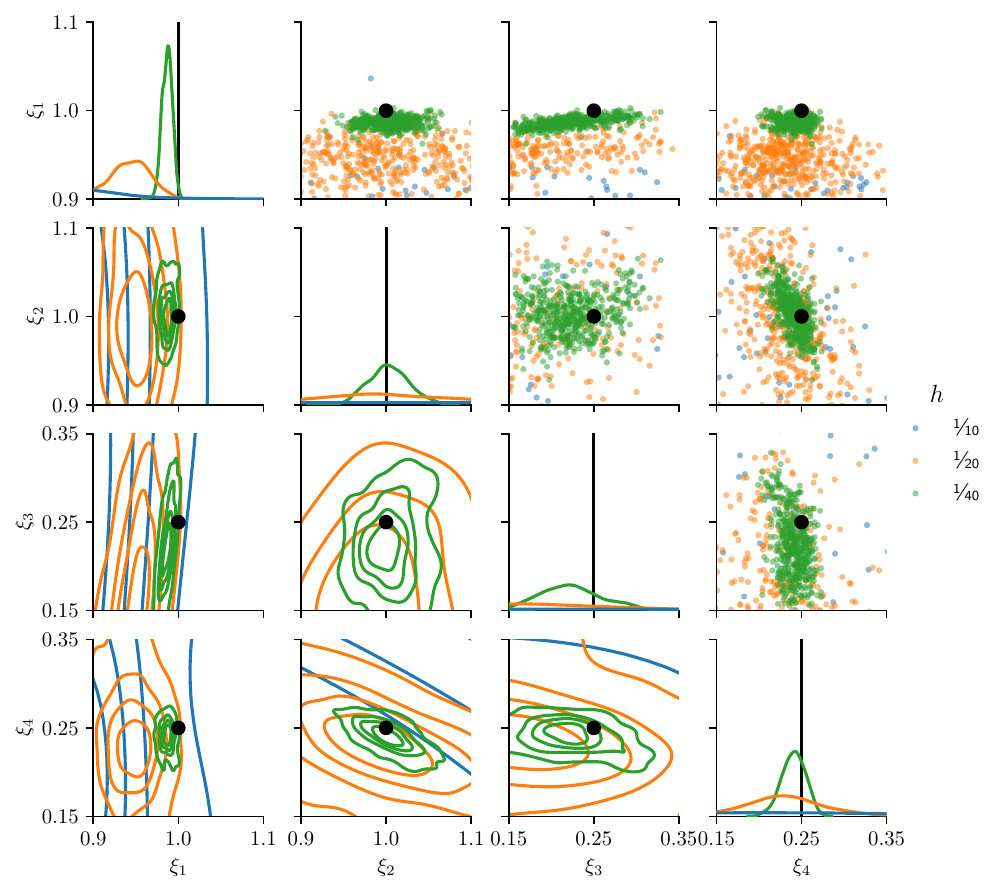}
        \caption{RM-FEM}
        \label{subfig:pairgrid-plot-rmfem}
    \end{subfigure}
    \caption{
        Pairgrid plot of the posterior $\posterior$ for both FEM and RM-FEM, using $\meshsize=\{\sfrac{1}{10}, \sfrac{1}{20}, \sfrac{1}{40}\}$.
        The true values of each parameter is indicated by a black line or dot.
    }
    \label{fig:pairgrid-plots}
\end{figure*}

In \cref{fig:kappa-plots}, the resulting posterior distributions $\posteriorkappa$ are shown, using both FEM and RM-FEM and for $\meshsize=\{\sfrac{1}{10}, \sfrac{1}{20}, \sfrac{1}{40}\}$.
To get a better insight in the fitted parameters, a pairgrid plot of $\posterior$ based on the same data is shown in \cref{fig:pairgrid-plots}.
The mean, standard deviation, as well as the error between the mean and ground truth are given in \cref{tbl:summary}.
Looking at the FEM results in \cref{subfig:pairgrid-plot-fem}, we can see that the width of the posterior distributions is not affected by the number of elements in the mesh $\nelem$.
This is not unexpected, since the shape of the posterior is determined only by the prior and the likelihood, neither of which accounts for the model misspecification due to discretization error.
This results in overconfidence of the model for $\meshsize=\sfrac{1}{10}$ and $\meshsize=\sfrac{1}{20}$, where the posterior probability density at the ground truth is essentially zero.
For $\meshsize=\sfrac{1}{40}$, the model no longer appears overconfident, since at this point the uncertainty due to discretization error has become sufficiently small compared to the uncertainty due to observation noise.
In \cref{tbl:summary}, we see that for $\meshsize=\sfrac{1}{40}$ the error in the posterior mean is less than 1 standard deviation for each $\param_i$.

Turning to the RM-FEM results in \cref{subfig:pairgrid-plot-rmfem}, we can see that the number of elements used in the model is affecting the width of the posterior distribution that is obtained.
Indeed, as $\nelem$ increases, the posterior mean moves in the direction of the ground truth, and the width of the posterior decreases accordingly.
However, comparing the width of the RM-FEM posterior in \cref{subfig:pairgrid-plot-rmfem} to the FEM posterior error in \cref{subfig:pairgrid-plot-fem}, the model now appears to be underconfident.
Taking $\probabilityof[\observations]{\param_2}$ as an example with $\meshsize=\sfrac{1}{40}$, we see in \cref{tbl:summary} that RM-FEM yields a 4-fold increase in the standard deviation compared to FEM (from \num{0.0055} to \num{0.0234}), despite the original FEM error of \num{0.035} being less than 1 standard deviation.
This increase in variance despite convergence having been reached can be seen as a form of underconfidence of the model.

More worryingly, the posterior mean is further away from the ground truth for the RM-FEM posterior than for the FEM posterior.
For $\meshsize=\sfrac{1}{10}$, the RM-FEM posterior mean falls outside the domain of the pairgrid plots for all $\xi_i$.
Looking at $\probabilityof[\observations]{\param_1}$ with $\meshsize=\sfrac{1}{20}$, the absolute error of the posterior mean has increased from \num{0.0038} for FEM to \num{0.0563} for RM-FEM.
Even for $\meshsize=\sfrac{1}{40}$, which is sufficient for an accurate FEM posterior, RM-FEM introduces an additional error to the posterior mean.

\begin{table*}[tbp]
    \centering
    \captionsetup{width=0.8\textwidth}
    \caption{Summary of the mean and standard deviation of the posterior $\posterior$ as well as the error of the posterior mean $|\refparams - \expectationof{\posterior}|$, for all FEM approaches investigated in this paper.}
    \label{tbl:summary}
    \begin{tabular}{llrrrrrrrrrrrr}
        \cmidrule[\heavyrulewidth]{3-14}
        \multicolumn{2}{c}{} & \multicolumn{3}{c}{1D/2D FEM} & \multicolumn{3}{c}{1D RM-FEM} & \multicolumn{3}{c}{1D RM-FEM, fix obs} & \multicolumn{3}{c}{2D RM-FEM} \\
        \cmidrule(lr){3-5} \cmidrule(lr){6-8} \cmidrule(lr){9-11} \cmidrule(lr){12-14}
        \multicolumn{1}{c}{}     & $\meshsize$ & mean & std & error & mean & std & error & mean & std & error & mean & std & error \\
        \midrule
        \multirow{3}{*}{$\param_1$} & $\sfrac{1}{10}$ & .\textcolor{gray}{}9848 & .\textcolor{gray}{00}13 & .\textcolor{gray}{0}152 & .\textcolor{gray}{}8207 & .\textcolor{gray}{0}900 & .\textcolor{gray}{}1793 & .\textcolor{gray}{}9805 & .\textcolor{gray}{00}52 & .\textcolor{gray}{0}195 & .\textcolor{gray}{}9064 & .\textcolor{gray}{0}428 & .\textcolor{gray}{0}936 \\
                                    & $\sfrac{1}{20}$ & .\textcolor{gray}{}9962 & .\textcolor{gray}{00}13 & .\textcolor{gray}{00}38 & .\textcolor{gray}{}9437 & .\textcolor{gray}{0}231 & .\textcolor{gray}{0}563 & .\textcolor{gray}{}9944 & .\textcolor{gray}{00}21 & .\textcolor{gray}{00}56 & .\textcolor{gray}{}9686 & .\textcolor{gray}{0}114 & .\textcolor{gray}{0}314 \\
                                    & $\sfrac{1}{40}$ & .\textcolor{gray}{}9991 & .\textcolor{gray}{00}13 & .\textcolor{gray}{000}9 & .\textcolor{gray}{}9869 & .\textcolor{gray}{00}61 & .\textcolor{gray}{0}131 & .\textcolor{gray}{}9986 & .\textcolor{gray}{00}14 & .\textcolor{gray}{00}14 & .\textcolor{gray}{}9926 & .\textcolor{gray}{00}33 & .\textcolor{gray}{00}74 \\
        \hline
        \multirow{3}{*}{$\param_2$} & $\sfrac{1}{10}$ & .\textcolor{gray}{}9479 & .\textcolor{gray}{00}57 & .\textcolor{gray}{0}521 & .\textcolor{gray}{}9820 & .\textcolor{gray}{}4145 & .\textcolor{gray}{0}180 & .\textcolor{gray}{}9329 & .\textcolor{gray}{0}179 & .\textcolor{gray}{0}671 &1.\textcolor{gray}{}0039 & .\textcolor{gray}{}1877 & .\textcolor{gray}{00}39 \\
                                    & $\sfrac{1}{20}$ & .\textcolor{gray}{}9871 & .\textcolor{gray}{00}54 & .\textcolor{gray}{0}129 & .\textcolor{gray}{}9947 & .\textcolor{gray}{0}932 & .\textcolor{gray}{00}53 & .\textcolor{gray}{}9807 & .\textcolor{gray}{00}66 & .\textcolor{gray}{0}193 & .\textcolor{gray}{}9887 & .\textcolor{gray}{0}457 & .\textcolor{gray}{0}113 \\
                                    & $\sfrac{1}{40}$ & .\textcolor{gray}{}9965 & .\textcolor{gray}{00}55 & .\textcolor{gray}{00}35 &1.\textcolor{gray}{}0035 & .\textcolor{gray}{0}234 & .\textcolor{gray}{00}35 & .\textcolor{gray}{}9946 & .\textcolor{gray}{00}57 & .\textcolor{gray}{00}54 & .\textcolor{gray}{}9983 & .\textcolor{gray}{0}123 & .\textcolor{gray}{00}17 \\
        \hline
        \multirow{3}{*}{$\param_3$} & $\sfrac{1}{10}$ & .\textcolor{gray}{}2349 & .\textcolor{gray}{00}89 & .\textcolor{gray}{0}151 & .\textcolor{gray}{}1913 & .\textcolor{gray}{}6465 & .\textcolor{gray}{0}587 & .\textcolor{gray}{}2264 & .\textcolor{gray}{0}243 & .\textcolor{gray}{0}236 & .\textcolor{gray}{0}845 & .\textcolor{gray}{}2961 & .\textcolor{gray}{}1655 \\
                                    & $\sfrac{1}{20}$ & .\textcolor{gray}{}2433 & .\textcolor{gray}{00}96 & .\textcolor{gray}{00}67 & .\textcolor{gray}{}1133 & .\textcolor{gray}{}1239 & .\textcolor{gray}{}1367 & .\textcolor{gray}{}2408 & .\textcolor{gray}{0}118 & .\textcolor{gray}{00}92 & .\textcolor{gray}{}1834 & .\textcolor{gray}{0}679 & .\textcolor{gray}{0}666 \\
                                    & $\sfrac{1}{40}$ & .\textcolor{gray}{}2453 & .\textcolor{gray}{00}85 & .\textcolor{gray}{00}47 & .\textcolor{gray}{}2221 & .\textcolor{gray}{0}367 & .\textcolor{gray}{0}279 & .\textcolor{gray}{}2432 & .\textcolor{gray}{00}96 & .\textcolor{gray}{00}68 & .\textcolor{gray}{}2352 & .\textcolor{gray}{0}207 & .\textcolor{gray}{0}148 \\
        \hline
        \multirow{3}{*}{$\param_4$} & $\sfrac{1}{10}$ & .\textcolor{gray}{}2230 & .\textcolor{gray}{00}46 & .\textcolor{gray}{0}270 & .\textcolor{gray}{}1433 & .\textcolor{gray}{}3002 & .\textcolor{gray}{}1067 & .\textcolor{gray}{}2175 & .\textcolor{gray}{0}213 & .\textcolor{gray}{0}325 & .\textcolor{gray}{}1713 & .\textcolor{gray}{}1247 & .\textcolor{gray}{0}787 \\
                                    & $\sfrac{1}{20}$ & .\textcolor{gray}{}2418 & .\textcolor{gray}{00}41 & .\textcolor{gray}{00}82 & .\textcolor{gray}{}2287 & .\textcolor{gray}{0}540 & .\textcolor{gray}{0}213 & .\textcolor{gray}{}2387 & .\textcolor{gray}{00}75 & .\textcolor{gray}{0}113 & .\textcolor{gray}{}2360 & .\textcolor{gray}{0}278 & .\textcolor{gray}{0}140 \\
                                    & $\sfrac{1}{40}$ & .\textcolor{gray}{}2474 & .\textcolor{gray}{00}40 & .\textcolor{gray}{00}26 & .\textcolor{gray}{}2405 & .\textcolor{gray}{0}137 & .\textcolor{gray}{00}95 & .\textcolor{gray}{}2470 & .\textcolor{gray}{00}42 & .\textcolor{gray}{00}30 & .\textcolor{gray}{}2448 & .\textcolor{gray}{00}77 & .\textcolor{gray}{00}52 \\
        \bottomrule
\end{tabular}
\end{table*}

\section{Interpolation error}
\label{sec:interpolation-error}
The explanation for the underconfidence of RM-FEM and the deterioration of its posterior mean has a local and a global component.
In this section, we will focus on the local component, interpolation error, by considering its relative contribution to the total discretization error.
Since our observations $\observations$ are point evaluations of the solution field $\solution\ofx$, we will consider the $L_2$ norm of the displacement error:
\begin{equation}
    \fielderrornorm = \left(\int_\domain \left(\refsolution\ofx - \femsolution\ofx \right)^2 \diff \domain\right)^\frac{1}{2}
\end{equation}
We will use a weighted norm for the nodal error of the finite element solution:
\begin{equation}
    \nodalerrornorm = \left( \sum_{i = 1}^{m} h_i \left(\refsolution_i - \femsolution_i\right)^2 \right)^\frac{1}{2}
\end{equation}
and define the relative nodal error $\nodalerrorrel$ and the relative interpolation error $\intererrorrel$ as:
\begin{align}
    \nodalerrorrel &= \frac{\nodalerrornorm}{\fielderrornorm} & \intererrorrel &= 1 - \nodalerrorrel
\end{align}

The nodal discretization error $\nodalerror$ is relatively small for FEM solutions to smooth elliptic 1D differential equations.
In fact, for a constant $\diffusion\ofx$, the finite element solution is exact in the nodes and the relative nodal error $\nodalerrorrel$ and interpolation error $\intererrorrel$ would be exactly \num{0} and \num{1}, respectively.
Although for the problem at hand, the nodal FEM solution is not exact, the relative interpolation error is still close to 1.

\begin{figure}
    \includegraphics[width=\columnwidth]{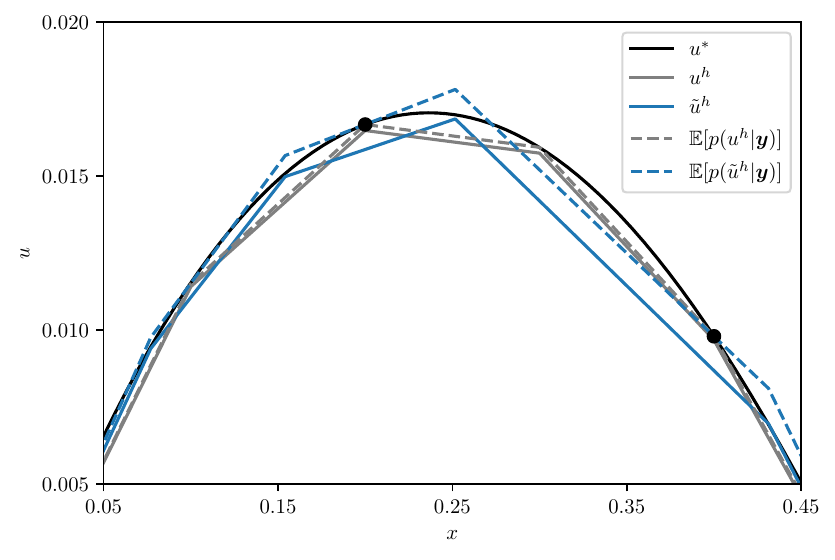}
    \caption{
        The FEM solution to \cref{eq:strong-form} on the reference mesh and a single perturbed mesh (in gray and blue, respectively).
        The posterior mean that is obtained when solving the Bayesian inference problem on either mesh is shown with a dashed line.
        The reference solution $\refsolution\ofx$ and observations $\observations$ are shown in black.
    }
    \label{fig:interpolation-plot}
\end{figure}

In \cref{fig:interpolation-plot}, the reference solution $\refsolution\ofx$ is shown in black as well as the FEM solution $\femsolution\ofx$ of the unperturbed mesh with $\meshsize=\sfrac{1}{10}$ and that of a single perturbed mesh $\pertfemsolution\ofx$ in gray.
For both FEM solves, the ground truth parameters $\refparams$ have been used.
Since the observation locations coincide with the unperturbed node locations, the error of the unperturbed FEM prediction in the observation locations is significantly smaller than the error of the perturbed FEM prediction.

We now solve the inverse problem on both meshes using the procedure described in \cref{subsec:mcmc}, and plot the resulting posterior mean in blue in \cref{fig:interpolation-plot}.
The smaller the distance between the FEM solution given the reference parameters $\femsolution\ofx|\refparams$ and the posterior (mean) FEM solution $\expectationof{\probabilityof[\observations]{\femsolution\ofx}}$, the smaller the error in the likelihood will be and therefore the smaller the effect of discretization error on the posterior distribution $\posterior$.
For this 1D problem, the relative nodal error $\nodalerrorrel$ is close to 0, and as a result the FEM solution is quite accurate, even for $\meshsize=\sfrac{1}{10}$.
This explains why the inverse solution using FEM shown in \cref{subfig:kappa-plot-fem-10,subfig:pairgrid-plot-fem} is very accurate, even for $\meshsize=\sfrac{1}{10}$.

For the perturbed mesh, we see that interpolation error significantly affects the accuracy of the prediction at the observation locations.
As a result, large mismatch arises between the perturbed FEM solution given the reference parameters $\pertfemsolution\ofx | \refparams$ and the posterior mean of the perturbed mesh $\expectationof{\probabilityof[\observations]{\pertfemsolution\ofx}}$.
Note that this effect will occur for any perturbed mesh:
the nodal locations of the unperturbed mesh are optimal in the sense that they minimize the FEM error at the observation locations, and consequently every perturbed mesh will be suboptimal in this regard.
This in turn causes a larger error in the likelihood, which translates to both the underconfidence and the deterioration of the posterior mean we observed in \cref{subfig:pairgrid-plot-rmfem}.

\subsection{Fixed observation nodes}
\label{subsec:fixed-observations}
Although worrisome, this issue is not necessarily insurmountable.
One potential approach to resolving it is to perturb only the nodes in the mesh that are not at one of the observation locations.
Although this would be problematic for $\meshsize=\sfrac{1}{5}$, in which case all nodes are observed and so none are perturbed, this approach could be feasible for finer meshes.
The smaller $\meshsize$ is, the smaller the effect of fixing the 4 observation nodes is expected to be.

\begin{figure*}
    \begin{subfigure}{0.5\textwidth}
        \includegraphics[width=\textwidth]{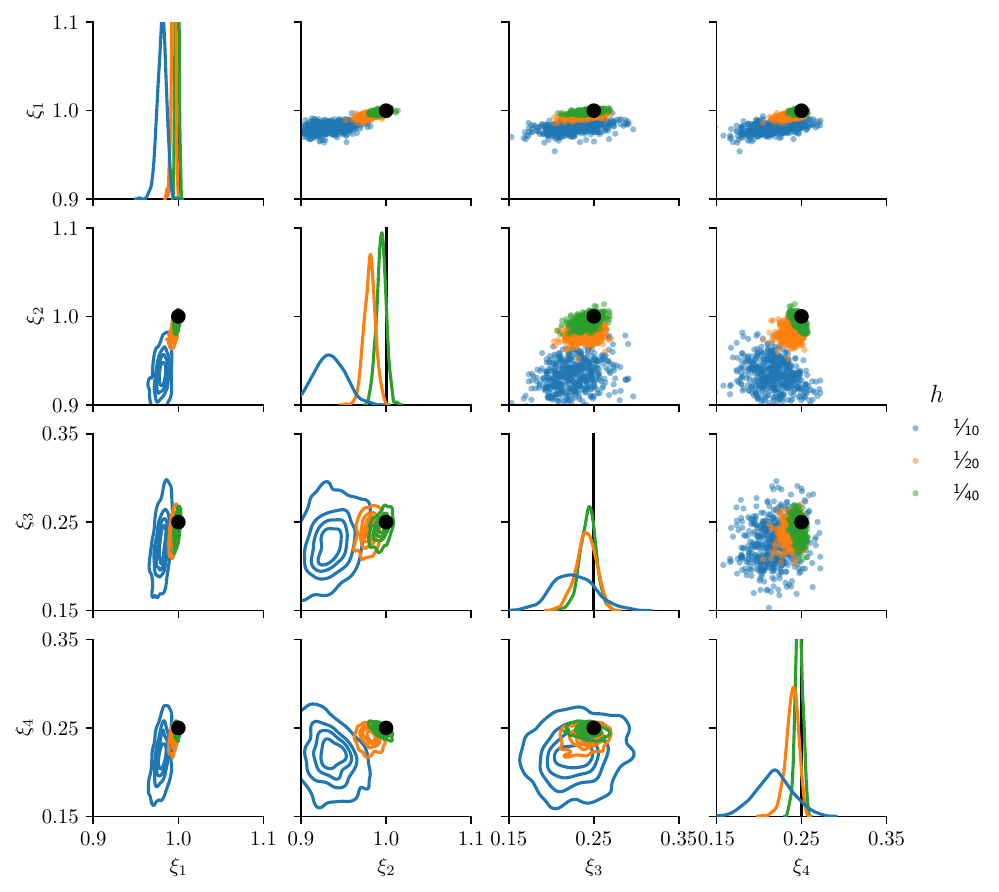}
        \caption{RM-FEM, fixed observation nodes}
        \label{subfig:pairgrid-plot-rmfem-omit}
    \end{subfigure}
    \hfill
    \begin{subfigure}{0.5\textwidth}
        \includegraphics[width=\textwidth]{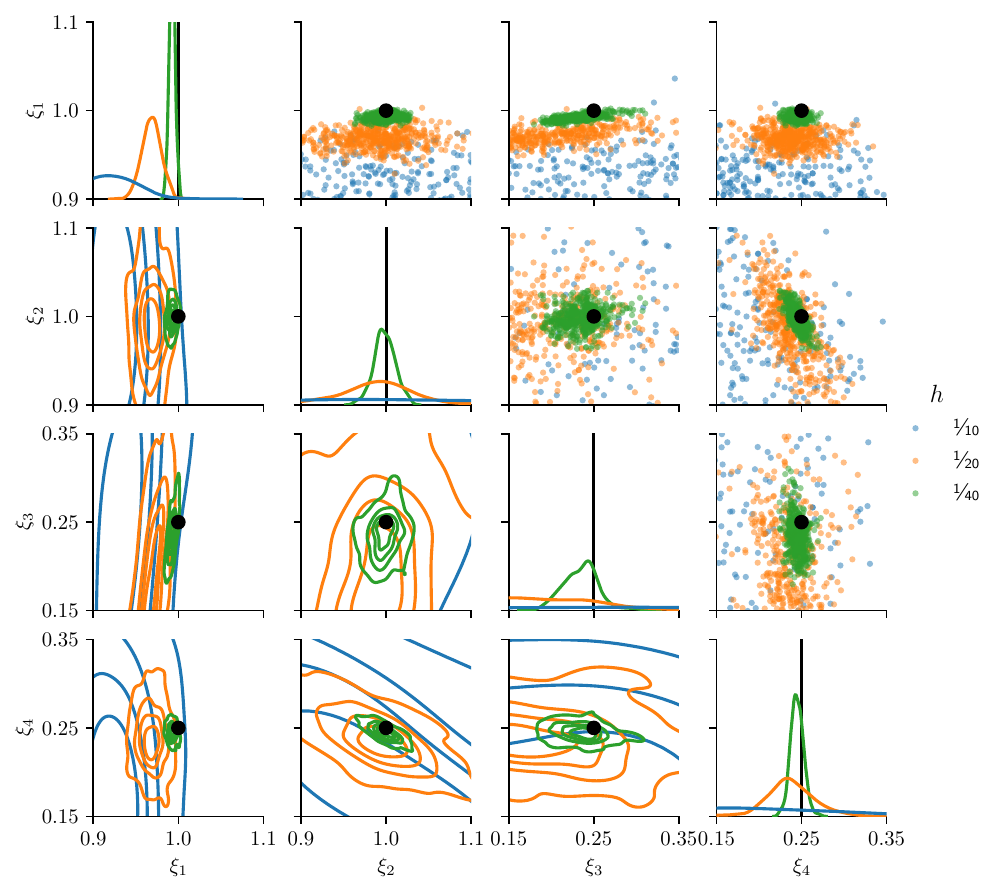}
        \caption{RM-FEM, 2D}
        \label{subfig:pairgrid-plot-rmfem-2d}
    \end{subfigure}
    \caption{
        Pairgrid plot of the posterior $\posterior$ for 1D RM-FEM with fixed observation nodes, and 2D RM-FEM, using $\meshsize=\{\sfrac{1}{10}, \sfrac{1}{20}, \sfrac{1}{40}\}$.
        The true parameter values are indicated by a black line or dot.
    }
    \label{fig:rmfem-fixes}
\end{figure*}

The resulting posterior $\posterior$ is shown in \cref{subfig:pairgrid-plot-rmfem-omit}, and its moments are summarized in \cref{tbl:summary} under ``1D RM-FEM, fix obs''.
We observe that fixing the observed nodes reduces the width of the RM-FEM posterior for all mesh densities $\meshsize$ and parameters $\param_i$, while maintaining the desired increase in confidence with mesh refinement.
Additionally, although a slight increase in the error is still visible compared to the FEM posterior mean, the severe deterioration of the RM-FEM posterior mean has been mitigated by fixing the observation nodes.
Furthermore, looking again at $\probabilityof[\observations]{\param_2}$ with $\meshsize=\sfrac{1}{40}$, we see hardly any difference in the standard deviation of the RM-FEM posterior with fixed observation nodes, compared to the FEM posterior.
This is a desirable property, since for this mesh density, the FEM error is small compared to the uncertainty due to observation noise.
We conclude that fixing the observed nodes helps reduce both the underconfidence of RM-FEM and the deterioration of its posterior mean.

\subsection{Generalizability}
\label{subsec:2D-problems}
We note that the conclusion of the previous section is quite specific to the 1D case with nodal observations.
If observations are made over a small region rather than a single point, the effect of interpolation error would be dampened.
Similarly, this relatively large interpolation error $\intererrorrel$ is typical for FEM solutions to elliptic PDEs in 1D, but not in 2D or 3D.
For more realistic problems, most of the discretization error is due to nodal error $\nodalerrorrel$ rather than interpolation error $\intererrorrel$.
Finally, higher-dimensional problems require a different perturbation scheme to guarantee a valid finite element discretization.
For example, in the case of quadrilateral elements, the maximal perturbation that still guarantees mesh validity is a disk centered on the nodes with a radius $r = \sfrac{1}{4} \sqrt{2} h$.
Considering only the $x$-component, we would obtain a semicircular distribution on the support $(-\sfrac{1}{4}\sqrt{2} h, \sfrac{1}{4}\sqrt{2} h)$ rather than a uniform distribution on the support $(-\sfrac{1}{2} h, \sfrac{1}{2} h)$.

Although all three of these factors are expected to reduce the relative effect of interpolation error in a more realistic setting, we will single out the final point (\ie a different perturbation scheme) to put it to the test.
We change to a 2D domain $\domain = (0, 1) \times (0, \sfrac{1}{10})$, with a Dirichlet boundary on the left and right $\boundarydir = \{0, 1\} \times (0, \sfrac{1}{10})$ and a Neumann boundary on the top and bottom $\boundaryneu = (0, 1) \times \{0, \sfrac{1}{10}\}$, such that $\boundarydir \cup \boundaryneu = \boundary$ and $\boundarydir \cap \boundaryneu = \emptyset$.
We adapt the PDE such that $\diffusion\ofvx$, $\forcing\ofvx$ and $\solution\ofvx$ are all constant with respect to the vertical coordinate $y$, and the reference solution $\refsolution$ as a function of $x$ is exactly the same as the solution to \cref{eq:strong-form}:
\begin{equation}
    \label{eq:strong-form-2d}
    \begin{aligned}
        - \divergence \left(\diffusion\ofvx \gradient \solution \ofvx \right) &= \forcing \ofvx & \forall \vx &\in \domain \\
        \solution\ofvx &= 0 & \forall \vx &\in \boundarydir \\
        \frac{\partial u}{\partial y}\ofvx &= 0 & \forall \vx &\in \boundaryneu
    \end{aligned}
\end{equation}
The coarsest mesh consists of 10 square elements with side length $\meshsize=\sfrac{1}{10}$ and linear shape functions.
Each finer mesh is obtained by splitting each element into 4 smaller square elements.
For our RM-FEM implementation, we perturb each node by choosing from a uniform distribution on a disk with radius $r = \sfrac{1}{4} \sqrt{2} h$, centered on the node location on the reference mesh.
The boundary nodes are perturbed in the same manner, but projected back onto the boundary $\boundary$ to ensure that the geometry of the domain is identical for each mesh.
We define the observation locations for this 2D problem to be $\obslocations = \{(\sfrac{i}{5}, \sfrac{1}{20})\}_{i=1}^4$, and use the same observation values $\observations$ as in the 1D case.
When solving the inverse problem with FEM, we obtain exactly the same posterior distribution as in the 1D case, shown in \cref{subfig:pairgrid-plot-fem}.

However, when using RM-FEM for this 2D case, we obtain the posterior shown in \cref{subfig:pairgrid-plot-rmfem-2d}, which does differ from the 1D RM-FEM posterior shown in \cref{subfig:pairgrid-plot-rmfem}.
Similarly to 1D RM-FEM with fixed observation nodes, we see in \cref{tbl:summary} that for this 2D problem, the width of the RM-FEM posterior has reduced significantly compared to 1D RM-FEM, and that the increase of confidence with mesh refinement is maintained.
Furthermore, the error of the 2D RM-FEM posterior mean is smaller than that of 1D RM-FEM for all parameters $\param_i$.
Despite the difference between 1D and 2D RM-FEM not being as great as the difference between 1D RM-FEM with and without fixed observation nodes, we can conclude that the 2D perturbation scheme reduces the severity of the deterioration of the posterior mean.
We expect that the deterioration of the posterior mean would vanish altogether for problems with regional observations (as opposed to point observations) and problems where the relative nodal error $\nodalerrorrel$ dwarfs the relative interpolation error $\intererrorrel$.

\section{Bias}
\label{sec:bias}
Although we have seen in the previous section that the adverse effects of interpolation error on the performance of RM-FEM can be dampened to some extent, some points remain unexplained.
In particular, it is still unclear why RM-FEM does not necessarily improve the accuracy of the parameter estimate provided by the posterior mean compared to FEM.
For instance, throughout all examples shown in \cref{fig:pairgrid-plots,fig:rmfem-fixes}, $\param_1$ can be seen to be consistently underestimated.
As we will see, the reason for this lies not in some local effect like the one we explored in the previous section, but rather lies in a global bias in the FEM estimate of the solution.

To explain this global effect, we turn back to the forward problem presented in \cref{subsec:forward-model} and consider the total energy $\refenergy$ in the system as well as its FEM approximation $\femenergy$:
\begin{align}
    \refenergy &= \int_\domain \solution\ofx \forcing\ofx \dx & \femenergy &= \T{\vecsolution} \vecforcing
\end{align}
FEM always underestimates this quantity, as a direct result of the fact that FEM searches the optimal solution in the subspace spanned by the shape functions.
This implicitly neglects the energy in the orthogonal space, and thus produces an estimate of the total energy $\refenergy$ that is too low.

\begin{figure}
    \centering
    \includegraphics[width=\columnwidth]{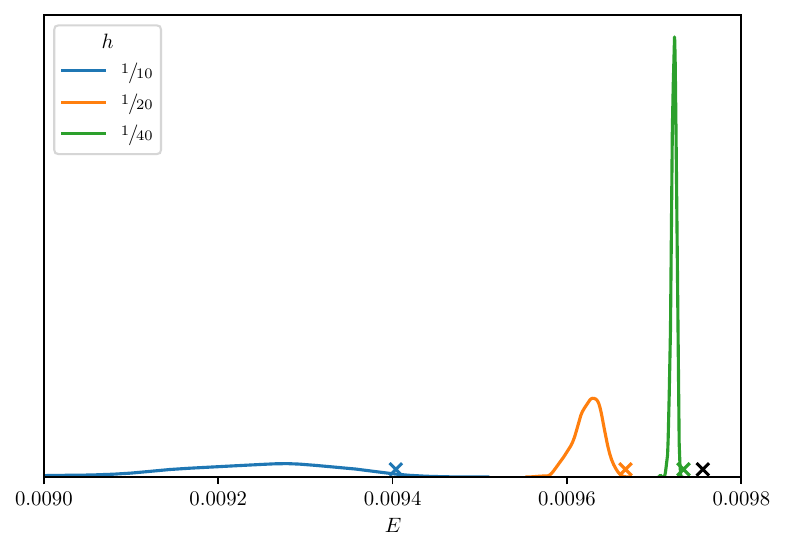}
    \caption{
        RM-FEM distributions of the total energy $\pertfemenergy$ for $\meshsize=\{\sfrac{1}{10}, \sfrac{1}{20}, \sfrac{1}{40}\}$.
        The total energy of the unperturbed meshes $\femenergy$ are indicated by colored crosses.
        The reference total energy $\refenergy$ is indicated by a black cross.
    }
    \label{fig:energy-norm-plot}
\end{figure}

We now compute $\femenergy$ using both FEM and RM-FEM for $\meshsize=\{\sfrac{1}{10}, \sfrac{1}{20}, \sfrac{1}{40}\}$, and show the results in \cref{fig:energy-norm-plot}.
As expected, the total energy $\refenergy$ is underestimated by FEM for $\meshsize=\sfrac{1}{10}$, and approaches the reference value as the mesh is refined.
We see this same underestimation occur when looking at the distributions produced by RM-FEM:
none of the perturbed meshes could ever produce a total energy estimate $\femenergy$ that is higher than or equal to the true total energy $\refenergy$.
This underestimation of the total energy typically results in FEM solutions that are biased toward 0, and this is compensated in the inverse setting by underestimating $\param_1$.

An additional observation to make about \cref{fig:energy-norm-plot} is that the energy estimate provided by FEM is close to the optimal energy of the RM-FEM distribution.
This is in line with common knowledge that smooth 1D differential equations like the one studied here are best discretized by a uniform mesh.
More broadly, the aim of discretization strategies is to choose a mesh such that $\femenergy$ is as close to $\refenergy$ as possible for a given number of nodes or elements.
A fundamental limitation of RM-FEM is that it pays little heed to the care that is typically put into mesh selection.
If the original FEM mesh was chosen arbitrarily, then the RM-FEM mean energy $\expectationof{\pertfemenergy}$ would on average coincide with the FEM estimate $\femenergy$.
However, since the FEM mesh is chosen close to optimally, we see a deterioration of the RM-FEM mean energy $\expectationof{\pertfemenergy}$ when comparing against the FEM estimate $\femenergy$.

We emphasize that these results do not contradict any of the conclusions drawn in \cite{abdulle_probabilistic_2021}.
The error estimators derived from the distribution over random meshes, as well as their application to adaptive mesh refinement are not implicated by any part of our study.
In fact, it may well be possible to derive an error estimator for $\refenergy$ from the RM-FEM distribution given in \cref{fig:energy-norm-plot}, since the test for a valid error estimator is whether it shares the same order of convergence as the true error.
From this point of view, whether $\refenergy$ falls inside or outside the $95\%$ confidence interval of the distribution is irrelevant.
However, when attempting to derive confidence bounds for a quantity of interest, or when applying RM-FEM in an inverse setting, caution is required.

We can think of FEM as a point estimator of the solution $\refsolution\ofx$ and of RM-FEM as a bootstrapping method\footnote{The term ``bootstrap'' is sometimes used specifically for case resampling. We are using the term here in a broader sense, \ie turning one dataset into $D$ datasets, fitting $D$ models, and leveraging their joint performance. This includes case resampling, but also encompasses other approaches (\eg the parametric bootstrap, the Bayesian bootstrap and residual resampling).}, whose moments provide useful information about the uncertainty of our point estimate.
The bootstrap can be a useful tool to estimate the variance of an estimator, and reduce the variance through bootstrap aggregation, but is often ill-equipped to estimate bias.
This carries over to RM-FEM, which is useful as an error estimator, but which produces unhelpful distributions for quantities of interest for which FEM is known to be biased.
The energy norm used in this section is a clear example of this and the bias we observed for $\param_1$ falls in the same category.
It is not obvious how one should correct such biases in RM-FEM.

\section{Conclusion}
\label{sec:conclusion}
In this work, we investigated the effects of interpolation error and bias on the performance of RM-FEM in an inverse setting.
In 1D, interpolation error can deteriorate the method's performance, yielding a posterior mean that is further away from the ground truth and an overly large variance, when compared to the posterior mean that is obtained using traditional FEM.
This effect can be dampened by excluding the nodes at the observation locations from the perturbation scheme, and is generally less pronounced in more practical settings.
A more pernicious drawback of RM-FEM is its quantification of discretization error for quantities of interest for which FEM is known to be biased, in which case the RM-FEM distribution is biased as well, possibly even more so.
Although RM-FEM makes a useful tool for error estimation and adaptive mesh refinement, we recommend caution when employing it in an inverse setting, particularly in cases where FEM produces biased estimates of the observed quantity of interest.

\bibliography{references}

\end{document}